# Density estimation with heteroscedastic error

AURORE DELAIGLE[1] and ALEXANDER MEISTER[2]

[1]*Department of Mathematics, University of Bristol, Bristol BS8 1TW, UK and Department of Mathematics and Statistics, University of Melbourne, VIC, 3010, Australia.*
*E-mail: aurore.delaigle@bristol.ac.uk*
[2]*Graduiertenkolleg 1100, Universität Ulm, Helmholtzstrasse 22, D-89081 Ulm, Germany.*
*E-mail: alexander.meister@uni-ulm.de*

It is common, in deconvolution problems, to assume that the measurement errors are identically distributed. In many real-life applications, however, this condition is not satisfied and the deconvolution estimators developed for homoscedastic errors become inconsistent. In this paper, we introduce a kernel estimator of a density in the case of heteroscedastic contamination. We establish consistency of the estimator and show that it achieves optimal rates of convergence under quite general conditions. We study the limits of application of the procedure in some extreme situations, where we show that, in some cases, our estimator is consistent, even when the scaling parameter of the error is unbounded. We suggest a modified estimator for the problem where the distribution of the errors is unknown, but replicated observations are available. Finally, an adaptive procedure for selecting the smoothing parameter is proposed and its finite-sample properties are investigated on simulated examples.

*Keywords:* bandwidth; density deconvolution; errors-in-variables; heteroscedastic contamination; inverse problems; plug-in

## 1. Introduction

We consider nonparametric estimation of a density from a sample contaminated by random error. This problem, which is called a *deconvolution problem*, arises very frequently in real data applications since, in practice, one often introduces non-negligible measurement errors while observing the data. The fields of application are various and include astronomy, biology, chemistry, economy and public health; see, for example, Merritt (1997) or the numerous examples described in Carroll *et al.* (2006).

In the conventional case, the observations are a sample of independent and identically distributed (i.i.d.) variables $Y_1, \ldots, Y_n$ generated by the model

$$Y_j = X_j + \varepsilon_j, \qquad X_j \sim f_X \quad \text{and} \quad \varepsilon_j \sim f_\varepsilon, \tag{1.1}$$

where the unknown density $f_X$ of $X_j$ is the quantity of interest, $\varepsilon_j$ are the error variables, independent of $X_j$, and $f_\varepsilon$ is known. In this context, Carroll and Hall (1988) and Stefanski and Carroll (1990) proposed the deconvolution kernel density estimator. Let $K$ be a square-integrable kernel function, $\omega_n > 0$ a smoothing parameter and, for all $t$, assume $f_\varepsilon^{\mathrm{ft}}(t) \neq 0$, where $g^{\mathrm{ft}}$ denotes the Fourier transform of a function $g$. The deconvolution kernel estimator is defined by

$$\tilde{f}_n(x) = \frac{1}{2\pi} \int \exp(-\mathrm{i}tx) K^{\mathrm{ft}}(t/\omega_n) \frac{1}{n} \sum_{j=1}^{n} \exp(\mathrm{i}tY_j)/f_\varepsilon^{\mathrm{ft}}(t) \, \mathrm{d}t; \qquad (1.2)$$

see, for example, Fan (1991a, b), Fan (1993) and Masry (1993) for theoretical properties. Recent contributions to density deconvolution include Zhang and Karunamuni (2000), Carroll and Hall (2004), van Es and Uh (2005), Hall and Qiu (2005) and Hall and Meister (2007).

In many applications of interest, the assumption of homoscedastic errors is too restrictive to be realistic. Bennett and Franklin (1954) describe an experiment where some students were asked to assess the iron content of substances. Here, clearly, the measurement process, and hence the error distribution, is subjective and differs among individuals. In some experiments, the error distribution depends on the type of individual under study (e.g., healthy or not, smoker or not, etc.) or on the measurement process. Here, as soon as the sample contains observations of different types, the errors are not identically distributed in the sample; see Fuller (1987) for an early consideration of this problem. Heteroscedasticity also arises when the sample is formed by collating data from different laboratories (see, e.g., National Research Council (1993)) or from different studies (meta-analysis), or when $r_i$ contaminated replications available for each individual $i$ are averaged to form a new sample of observations – a procedure often used in practice, because it reduces the scale of error.

In Section 2, we formally introduce the heteroscedastic error model and propose a deconvolution kernel estimator of the density $f_X$ that accounts for heteroscedastic errors. We establish $L_2$-consistency of the estimator, obtain its rates of convergence and prove that these are optimal. In Section 3, we study two important aspects of heteroscedastic contamination. We first consider the problem where different numbers of replicates are observed for each random variable $X_j$. We show that, in the case of normal contamination, averaging the replicates and then using the procedure derived in Section 2 leads to optimal convergence rates. Next, we discuss limiting cases of heteroscedastic errors with unbounded scaling parameters and give an equivalent criterion for the existence of a consistent estimator. Section 4 discusses some situations where the error distributions are unknown, but either replicated observations are available or more restrictive conditions on $f_X$ are assumed. We study finite-sample properties of our estimator in Section 5. We develop a data-driven bandwidth selector and give some numerical simulations. All proofs are deferred to Section 6.

# 2. Estimation procedure and asymptotics

## 2.1. The estimator

We generalize model (1.1) to allow heteroscedastic contamination, leading to the model

$$Y_j = X_j + \varepsilon_j, \qquad X_j \sim f_X \quad \text{and} \quad \varepsilon_j \sim f_{\varepsilon_j}. \tag{2.1}$$

Now, each $\varepsilon_j$ has its own density $f_{\varepsilon_j}$, which may depend on both the observation number $j$ and the sample size $n$. In this setting, where (1.2) can no longer be used, the estimator we propose is defined by

$$\widehat{f}_n(x) = \frac{1}{2\pi} \int \exp(-\mathrm{i}tx) K^{\mathrm{ft}}(t/\omega_n) \Psi_n(t) \, \mathrm{d}t \tag{2.2}$$

with

$$\Psi_n(t) = \sum_{j=1}^{n} f_{\varepsilon_j}^{\mathrm{ft}}(-t) \exp(\mathrm{i}tY_j) \bigg/ \left( \sum_{k=1}^{n} |f_{\varepsilon_k}^{\mathrm{ft}}(t)|^2 \right).$$

This estimator is well defined if we assume the following.

**Condition A.**

$$\text{There exists some } j \text{ such that } |f_{\varepsilon_j}^{\mathrm{ft}}(t)| \neq 0 \text{ for all } t \in \mathbb{R}, \tag{A.1}$$

$$K^{\mathrm{ft}}(t) \text{ is bounded, continuous at } t = 0 \text{ and } K^{\mathrm{ft}}(0) = 1, \tag{A.2}$$

$$f_{\varepsilon_j}^{\mathrm{ft}}(t\omega_n) K^{\mathrm{ft}}(t) \bigg/ \left( \sum_{k=1}^{n} |f_{\varepsilon_k}^{\mathrm{ft}}(t\omega_n)|^2 \right) \in L_2(\mathbb{R}) \qquad \text{for } j = 1, \dots, n. \tag{A.3}$$

These conditions are standard in deconvolution problems. In particular, in order to satisfy (A.3), it is rather common to choose kernels that have a compactly supported Fourier transform $K^{\mathrm{ft}}$. Such kernels are supported on the whole real line, examples being the sinc kernel $K_1(x) = \sin x/(\pi x)$ and the kernel $K_2(x) = 48(\cos x)(1 - 15x^{-2})/(\pi x^4) - 144(\sin x)(2 - 5x^{-2})/(\pi x^5)$, which have respective characteristic functions $K_1^{\mathrm{ft}}(t) = 1_{[-1,1]}(t)$, the indicator function of the interval $[-1, 1]$, and $K_2^{\mathrm{ft}}(t) = (1 - t^2)^3 1_{[-1,1]}(t)$.

An alternative estimator that can perhaps be seen as a more natural generalization of (1.2) is the estimator obtained when using $n^{-1} \sum_{j=1}^{n} \exp(\mathrm{i}tY_j) \{f_{\varepsilon_j}^{\mathrm{ft}}(t)\}^{-1}$ instead of $\Psi_n(t)$. A quick inspection of its properties, however, shows that this estimator suffers from the convergence rates of the least favorable error $\varepsilon_j$ and is therefore not acceptable. Another estimator of $f_X$, $\widehat{f}_{n,2}(x)$, can be defined if we replace $\Psi_n(t)$ by $\Phi_n(t) = \sum_{j=1}^{n} \exp(\mathrm{i}tY_j)/(\sum_{k=1}^{n} f_{\varepsilon_k}^{\mathrm{ft}}(t))$. As an advantage, applying this estimator requires only knowledge of the set $\{f_{\varepsilon_1}, \dots, f_{\varepsilon_n}\}$, but not the information about which observation is corrupted by which of the error densities. However, it is less attractive in some cases of non-symmetric $f_{\varepsilon_k}$ as, then, there is no guarantee that the denominator in $\Phi_n(t)$ does not vanish, although each $f_{\varepsilon_k}^{\mathrm{ft}}$ is assumed to have no zeros. Also, the mean integrated squared error of (2.2) is smaller than that of $\widehat{f}_{n,2}$ and, therefore, for the most part, we will focus our consideration on (2.2).

## 2.2. Asymptotic properties

We study asymptotic properties of our estimator by examining its mean integrated squared error (MISE), defined by $\mathrm{MISE}_n(f_X) = \mathrm{E}\|\widehat{f}_n - f_X\|^2_{L_2(\mathbb{R})}$. The usual bias-variance decomposition and the use of Parseval's identity lead to the following result.

**Lemma 2.1.** *Under Condition A, if $f_X \in L_2(\mathbb{R})$, then the estimator (2.2) satisfies*

$$
\begin{aligned}
\mathrm{MISE}_n(f_X) = {} & \frac{1}{2\pi} \int |f_X^{\mathrm{ft}}(t)|^2 |K^{\mathrm{ft}}(t/\omega_n) - 1|^2 \, \mathrm{d}t \\
& + \frac{1}{2\pi} \int |K^{\mathrm{ft}}(t/\omega_n)|^2 \left( \sum_{k=1}^{n} |f_{\varepsilon_k}^{\mathrm{ft}}(t)|^2 \right)^{-1} \mathrm{d}t \\
& - \frac{1}{2\pi} \int \left( \sum_{j=1}^{n} |f_{\varepsilon_j}^{\mathrm{ft}}(t)|^2 \right)^{-2} \left( \sum_{k=1}^{n} |f_{\varepsilon_k}^{\mathrm{ft}}(t)|^4 \right) |f_X^{\mathrm{ft}}(t)|^2 |K^{\mathrm{ft}}(t/\omega_n)|^2 \, \mathrm{d}t.
\end{aligned}
\tag{2.3}
$$

From the above lemma, we will be able to derive the rates of convergence of our estimator and prove their optimality in $\mathcal{F}_{\beta,C}$, the class of densities uniformly bounded relative to their Sobolev ($\beta$-)norm, that is, that satisfy

$$
\int |f_X^{\mathrm{ft}}(t)|^2 (1+t^2)^\beta \, \mathrm{d}t \le C.
\tag{2.4}
$$

Throughout, we assume $\beta > 1/2$, which ensures, for example, continuity of $f_X$. We also assume that the kernel $K$ satisfies the following condition, which is fulfilled by, for example, the sinc kernel $K_1$ (for any $\beta > 1/2$).

**Condition B.** $|K^{\mathrm{ft}}(t)| \le 1$ *for all $t$, $K^{\mathrm{ft}}$ is supported on $[-1, 1]$ and $|K^{\mathrm{ft}}(t) - 1| = o(|t|^\beta)$ with $\beta$ as in (2.4).*

Finally, we need some regularity assumptions on the error densities $f_{\varepsilon_j}$: we assume the existence of $\alpha, C > 0$ and the existence of some positive monotone decreasing functions $\overline{\varphi}_{j,n}(t)$ and $\underline{\varphi}_{j,n}(t)$ for $t > 0$ such that the following condition holds.

**Condition C.**

$$P(|\varepsilon_j| \leq \alpha) \geq C, \qquad \forall j, n, \tag{C.1}$$

$$|f^{\mathrm{ft}}_{\varepsilon_j}(t)| \geq \underline{\varphi}_{j,n}(T), \qquad \forall |t| \leq T, \tag{C.2}$$

$$\underline{\varphi}_{j,n}(t) \leq |f^{\mathrm{ft}}_{\varepsilon_j}(t)| \leq \overline{\varphi}_{j,n}(t), \qquad \forall t > T, \tag{C.3}$$

$$|f^{\mathrm{ft}}_{\varepsilon_j}{}'(t)| \leq \overline{\varphi}_{j,n}(t), \qquad \forall t > T, \tag{C.4}$$

$$\underline{\varphi}_{j,n}(t) \geq c_1 \cdot \overline{\varphi}_{j,n}(c_2 t), \qquad \forall t > 0, \tag{C.5}$$

*for some $T \geq 0$, $c_1 > 0$ and $c_2 \geq 1$ which are independent of $j$ and $n$. Note that condition (C.1) prevents $f_{\varepsilon_j}$ from spreading too intensively, while the other conditions represent a weak version of monotonicity for $|f^{\mathrm{ft}}_{\varepsilon_j}|$. In particular, the so-called ordinary smooth densities $f_{U_j}$, in the terminology of Fan (1991a, b), satisfy $\underline{\varphi}_{j,n}(t) = C_1 |t|^{-\nu}$ and $\overline{\varphi}_{j,n}(t) = C_2 |t|^{-\nu}$ with $C_2 > C_1 > 0$ and $\nu > 0$, and the supersmooth densities satisfy $\underline{\varphi}_{j,n}(t) = C_1 |t|^{\rho_1} \exp(-c|t|^\gamma)$ and $\overline{\varphi}_{j,n}(t) = C_2 |t|^{\rho_2} \exp(-c|t|^\gamma)$ with $C_2 \geq C_1 > 0$, $c > 0$, $\gamma > 0$ and $\rho_2 \geq \rho_1 \geq 0$.*

Under these conditions, we are ready to establish the rates of convergence of our estimator; the following theorem shows that, if the bandwidth is chosen appropriately, then our estimator achieves optimal rates.

**Theorem 2.1.** *Under Conditions A–C, assume the existence of a sequence $m_n \uparrow \infty$ such that, for some $C_2 \geq C_1 > 0$, $\beta > 1/2$,*

$$C_1 m_n^{1+2\beta} \leq \sum_{j=1}^n |\overline{\varphi}_{j,n}(m_n)|^2 \leq C_2 m_n^{1+2\beta} \tag{2.5}$$

*holds for all $n$. Then,*

(a) *when selecting $\omega_n = c_2^{-1} m_n$ (with $c_2$ defined in (C.5)), the estimator (2.2) fulfills*

$$\sup_{f_X \in \mathcal{F}_{\beta,C}} \mathrm{MISE}_n(f_X) = O(m_n^{-2\beta});$$

(b) *for an arbitrary estimator based on $Y_1, \ldots, Y_n$ and $C$ in (2.4) large enough, we have*

$$\sup_{f_X \in \mathcal{F}_{\beta,C}} \mathrm{MISE}_n(f_X) \geq const. \cdot m_n^{-2\beta}.$$

A more precise asymptotic description of the MISE, which we denote by AMISE, can be obtained under additional assumptions, by using a Taylor expansion of the bias term. Such an asymptotic expression is useful for deriving a data-driven bandwidth (see Section 5.1). Assume the following.

**Condition D.**

$$\omega_n \to \infty \text{ and } n/\omega_n \to \infty \text{ as } n \to \infty, \tag{D.1}$$

$$K \text{ is such that } \int |y^k K(y)| \, \mathrm{d}y < \infty \text{ and is of order } k, \tag{D.2}$$

$$f_X \text{ is } k+1 \text{ times differentiable, } \sup_{j=0,\ldots,k+1} \|f_X^{(j)}\|_\infty < \infty \text{ and } f_X^{(k)} \in L_2(\mathbb{R}), \tag{D.3}$$

*where a kth-order kernel is a kernel that satisfies* $\mu_{K,j} \equiv \int x^j K(x) \, \mathrm{d}x = 1_{\{j=0\}}$ *for* $j = 0,\ldots,k-1$ *and* $\mu_{K,k} = c$, *with* $c \neq 0$ *some finite constant.*

The AMISE is described in the next lemma, where we use the standard notation $h = \omega_n^{-1}$ for the bandwidth in order to highlight the usual bias-variance trade-off.

**Lemma 2.2.** *Under Conditions* A *and* D, *the estimator* (2.2) *satisfies* $\mathrm{MISE}_n(f_X) = \mathrm{AMISE}_n(f_X) - R_n + o(h^{2k})$, *where*

$$\mathrm{AMISE}_n(f_X) = \frac{h^{2k}\mu_{K,k}^2}{(k!)^2} \int (f_X^{(k)}(x))^2 \, \mathrm{d}x + \frac{1}{2\pi h} \int |K^{\mathrm{ft}}(t)|^2 \left(\sum_{k=1}^n |f_{\varepsilon_k}^{\mathrm{ft}}(t/h)|^2\right)^{-1} \mathrm{d}t \tag{2.6}$$

*and* $R_n = (2\pi)^{-1} \int (\sum_{j=1}^n |f_{\varepsilon_j}^{\mathrm{ft}}(t)|^2)^{-2} (\sum_{k=1}^n |f_{\varepsilon_k}^{\mathrm{ft}}(t)|^4) |f_X^{\mathrm{ft}}(t)|^2 |K^{\mathrm{ft}}(t/\omega_n)|^2 \, \mathrm{d}t.$

It can be shown that under mild conditions (e.g., Condition C), the term $R_n$ is negligible compared to the AMISE.

# 3. A few interesting results in limiting cases

This section is dedicated to studying a few interesting results obtained when considering limiting cases of model (2.1). We consider two extreme and opposite situations – error scales tending to zero or tending to infinity – and see how well the estimator behaves in these cases.

## 3.1. Averaging replicated observations

*Context.* Consider the rather frequent situation where the errors are homoscedastic and, for some individuals, replicated observations are available. The observations are of the form

$$Y_{j,k} = X_j + \varepsilon_{j,k}, \qquad j \in \{1,\ldots,n\}, k \in \{1,\ldots,r_{j,n}\}, \tag{3.1}$$

where $\varepsilon_{j,k} \sim f_\varepsilon$. When such data are available, it is rather common to work with the averaged observations $\overline{Y}_j = r_{j,n}^{-1} \sum_{k=1}^{r_{j,n}} Y_{j,k}$. Indeed, although, in (asymptotic) theory, using the averaged sample is not always advantageous – in some cases (ordinary smooth), the averaged errors become smoother and thus imply a slower rate of convergence – in finite samples, the variance reduction induced by the averaging process can lead to significant improvement of performance of the estimator; see Delaigle ([2008](#)). In this context, we apply our estimator ([2.2](#)) to the sample $\overline{Y}_j = X_j + \overline{\varepsilon}_j$, where, since $r_{j,n}$ may differ among individuals, the errors $\overline{\varepsilon}_j := r_{j,n}^{-1} \sum_{k=1}^{r_{j,n}} \varepsilon_{j,k}$ are heteroscedastic. Below, we denote the density of $\overline{\varepsilon}_j$ by $f_{\varepsilon_j}$.

*The normal case.* In many real data applications, it is reasonable to assume that the error is normally distributed, that is, $f_\varepsilon = N(\mu, \sigma^2)$ and $f_{\varepsilon_j} = N(\mu, \sigma_{j,n}^2)$ with $\sigma_{j,n}^2 = \sigma^2/r_{j,n}$ and $f_{\varepsilon_j}^{\mathrm{ft}}(t) = f_\varepsilon^{\mathrm{ft}}(t/r_{j,n})$. First, we show that in this case, there is no loss of information when using the averaged sample to estimate $f_X$.

**Theorem 3.1.** *Suppose $f_\varepsilon = N(\mu, \sigma^2)$ in the model ([3.1](#)). Then, the sample $\overline{Y}_1, \ldots, \overline{Y}_n$ is sufficient for $f_X$.*

It is clear that each $f_{\varepsilon_j}$ satisfies Condition [C](#); Conditions [A](#), [B](#) and ([2.5](#)) hold by appropriate selection of $K$ and $\omega_n$. Hence, Theorem [2.1](#) ensures rate optimality of our estimator ([2.2](#)) applied to the averaged data. It is not hard to prove that for $r_{j,n}$ fixed, the convergence rates of $\widehat{f}_n$ (when using the sample of averages, rather than the original sample) remain unchanged, but the constants improve (hence, the estimator behaves better with averaged data).

To gain more intuition about the amount of improvement one can get when using averaged data, consider the rather extreme situation where, as the sample size increases, more and more replicated data become available. The result below then shows that the usual logarithmic rates of convergence of the normal case can even become algebraic (see also Hesse ([1996](#)) for a related problem in the partial contamination context).

**Theorem 3.2.** *Under the conditions of Theorem [3.1](#), one is able to obtain algebraic rates for the supremum of the MISE taken over $f_X \in \mathcal{F}_{\beta,C}$, $\beta > 1/2$, if and only if there are some $\alpha > 0$, $\gamma > 0$, $c > 0$, $\delta > 0$ such that*

$$\# J_{n,\gamma,\alpha} \geq c \cdot n^\delta \qquad \text{for all } n, \tag{3.2}$$

*where we define $J_{n,\gamma,\alpha} := \{j \in \{1, \ldots, n\} : \sigma_{j,n}^2 < \gamma \cdot n^{-\alpha} \ln n\}$.*

For example, we easily verify ([3.2](#)) in the case $r_{j,n} \sim j^{\alpha_1} n^{\alpha_2}$ with $\alpha_1, \alpha_2 \geq 0$ and $\alpha_1 + \alpha_2 > 0$. Quite surprisingly, we notice the occurrence of algebraic rates in that case without the need for the total number of original data $N = \sum_{j=1}^n r_{j,n}$ to increase exponentially fast with increases of $n$. Here, $N$ increases only at a polynomial rate with $n$.

### 3.2. A case of unbounded scaling parameters

Whereas Theorem [3.2](#) focused on the behavior of our estimator in an extreme case where the error scale tends to zero, we now consider an opposite extreme situation where the scaling parameters are unbounded. We study this problem in the particular case where the $f_{\varepsilon_j}$ are symmetric and have the Fourier transform $f_{\varepsilon_j}^{\text{ft}}(t) = \exp(-\sigma_{j,n}^{\gamma}|t|^{\gamma}/2)$ with $\gamma \geq 1$ and some unbounded scaling parameters $\sigma_{j,n} > 0$. Examples of such densities are Cauchy densities for $\gamma = 1$ and centered normal densities for $\gamma = 2$, where $\sigma_{j,n}$ are scaling parameters. In this case, (C.1) is not satisfied and Theorem 2.1 cannot be applied. The next theorem shows the somewhat surprising result that if the unbounded sequence $(\sigma_{j,n})_{j,n}$ does not converge too rapidly to infinity, then the estimator remains consistent.

**Theorem 3.3.** (a) *With a suitable choice of $\omega_n$ and $K$ so that $K^{\text{ft}}$ is compactly supported and Condition A is satisfied, estimator (2.2) is consistent for $f_X$ without any smoothness assumptions on $f_X$ if, for any $\omega > 0$, we have*

$$\sum_{j=1}^{n} \exp(-\sigma_{j,n}^{\gamma}\omega^{\gamma}) \xrightarrow{n \to \infty} \infty. \tag{3.3}$$

(b) *If (3.3) is not valid, then there is no consistent estimator for $f_X \in \mathcal{F}_{\beta,C}$ with arbitrary $\beta > 1/2$ and $C$ large enough.*

This theorem also shows that the estimator (2.2) achieves consistency whenever consistent estimation is theoretically possible, for $\beta > 1/2$ and $C$ large enough. Examples of unbounded sequences that satisfy equation (3.3) are $\sigma_{j,n}^{\gamma} \leq o_n \cdot \log n$ and $\sigma_{j,n}^{\gamma} \leq o_j \cdot \log j$, where $o_n$ is an arbitrary sequence tending to zero.

# 4. The case of unknown error densities

Most papers dealing with deconvolution problems assume that the error densities are perfectly known as, otherwise, the target density is not identifiable in the standard models. However, since the error density is unknown in many practical situations, this classical condition is relaxed in some recent papers. As a payback, those models require either the availability of additional direct data from the error distribution (Diggle and Hall (1993), Neumann (1997)) or replicated measurements (Horowitz and Markatou (1996), Schennach (2004), Delaigle, Hall and Müller (2007), Delaigle, Hall and Meister (2008)) or more restrictive conditions on the target density (Butucea and Matias (2005), Meister (2006, 2007)).

In the heteroscedastic framework, the replicated measurement approach is of particular practical importance. In the context of Section 3.1, for example, that is, replicated measurement under normal contamination, and where the mean $\mu = 0$, but the variance $\sigma^2$ is unknown, $\sigma^2$ is estimable by $\widehat{\sigma^2} = (2N)^{-1} \sum_{(j,k_1,k_2) \in \mathcal{S}} (Y_{j,k_1} - Y_{j,k_2})^2$, where $\mathcal{S} = \{(j, k_1, k_2) \text{ such that } 1 \leq j \leq n, \ 1 \leq k_1 < k_2 \leq r_{j,n}\}$ and $N = \#\mathcal{S}$. The estimated variance $\widehat{\sigma^2}$ may replace $\sigma^2$ in the estimator (2.2) and it can be shown that this does not alter the convergence rates of Theorem 2.1 for sufficiently smooth $f_X$. This parametric procedure of error estimation is fairly standard in homoscedastic deconvolution because the possibility of obtaining replicated measurements is usually quite realistic; see, for example, Carroll, Eltinge and Ruppert (1993), Stefanski and Bay (1996), Carroll *et al.* (2004) and the references therein.

More surprisingly, in the most general case of our much less standard setting, where all error distributions are allowed to be different and no parametric shape is assumed for their densities, we are still able to use those replicates to consistently estimate the density $f_X$ under certain smoothness conditions. Indeed, if each observation is replicated at least once, $f_X$ can be consistently estimated by $\widehat{f}_{n,2}$ introduced in Section 2, where $\Phi_n(t)$ is replaced by the nonparametric estimator

$$\widehat{\Phi}_n(t) = \sum_{(j,k_1,k_2)\in\mathcal{S}} \exp(\mathrm{it}\{Y_{j,k_1}+Y_{j,k_2}\}/2) \Big/ \left( \left| \sum_{(k,k_1,k_2)\in\mathcal{S}} \sum_{(j,k_1,k_2)\in\mathcal{S}} \exp\{\mathrm{it}(Y_{j,k_1}-Y_{j,k_2})\} \right| + \rho \right),$$

with $\rho > 0$ a ridge parameter introduced to avoid division by zero and $\mathcal{S}$ as above. For symmetric error densities with non-vanishing Fourier transforms and appropriate selection of $h$ and $\rho$, consistency remains valid, even if the replicates of the same $X_j$ have different error distributions. We note, however, that in this very general case, the convergence rates of Theorem 2.1 cannot be maintained when the errors are ordinary smooth.

In the homoscedastic case, if the error density is known up to a scaling parameter, it is sometimes possible to estimate both that parameter and the target density $f_X$ without replicates. However, this can only be done by imposing more restrictive conditions on $f_X$ since a specific lower bound on $f_X^{\mathrm{ft}}$ has to be assumed (see Butucea and Matias (2005) or Meister (2006)). Under some circumstances, such methods can be extended to the heteroscedastic problem. For example, suppose we can assume that $f_X$ is symmetric and satisfies $|f_X^{\mathrm{ft}}(t)| \geq c/(1+|t|^{\beta+1/2})$ for all $t\in\mathbb{R}$ and some known $\beta>0$ and $c>0$, and that each error $\varepsilon_j$ is $N(0,\sigma_j^2)$, where $\sigma_j^2 = a(1+j/n)$, that is, the error variances follow a linear model with an unknown parameter $a$, say, in $[1,2]$. Note that $\varphi(a,t) \equiv n^{-1}\sum_{j=1}^n \exp(-\sigma_j^2 t^2/2) f_X^{\mathrm{ft}}(t)$ is $n^{-1}$-consistently estimable by the maximum of zero and the real part of the empirical characteristic function of the data for any $t$. Define known upper and lower bounds on $\varphi(a,t)$ by $\overline{\varphi}(a,t) = n^{-1}\sum_{j=1}^n \exp(-\sigma_j^2 t^2/2)$ and $\underline{\varphi}(a,t) = n^{-1}\sum_{j=1}^n \exp(-\sigma_j^2 t^2/2) c/(1+|t|^{\beta+1/2})$, respectively. We notice that for any $a > a'$, we have $\overline{\varphi}(a,t) < \underline{\varphi}(a',t)$ for $t$ sufficiently large. Introducing an equidistant partition of the interval $[1,2]$, where $a_j = 1+j/m$, $j = 1,\dots,m$, are the grid points, we fix $t$ large enough so that $\overline{\varphi}(a_{j-1},t) > \underline{\varphi}(a_{j-1},t) > \overline{\varphi}(a_j,t) > \underline{\varphi}(a_j,t) > \overline{\varphi}(a_{j+1},t) > \underline{\varphi}(a_{j+1},t)$. If, for some $j$, the empirically accessible function $\varphi(a_j,t)$ lies between $\overline{\varphi}(a_j,t)$ and $\overline{\varphi}(a_{j+1},t)$ we have $a \in [a_{j-1},a_{j+1}]$ as $\varphi(a,t)$ decreases monotonically in $a$. Then, by setting $m \to \infty$ at an appropriate order in $n$, we are able to estimate $a$; we may then insert its empirical counterpart $\widehat{a}$ into the estimator (2.2). Although those identification methods are very interesting, the framework of the current paper does not allow a more comprehensive study of this problem. However, we have learned that it is sometimes possible to extend the basic ideas of Butucea and Matias (2005) and Meister (2006) to the heteroscedastic setting.

# 5. Finite-sample performance

## 5.1. Data-driven bandwidth selection

We define the optimal bandwidth as the one that minimizes the MISE and estimate this bandwidth by a plug-in method similar to Delaigle and Gijbels (2004). We follow along the lines of their two-stage procedure and only explain the differences with their estimator, for a $k$th-order kernel. We select the bandwidth that minimizes the estimator of the AMISE in (2.6), obtained by replacing the unknown quantity $\theta_k = \int \{f_X^{(k)}\}^2$ by $\widehat{\theta}_k = \int \{\widehat{f}_n^{(k)}\}^2$, where, for $r$ any positive integer, $\widehat{f}_n^{(r)}(x) = (2\pi)^{-1} \int (-\mathrm{i}t)^r \exp(-\mathrm{i}tx) K^{\mathrm{ft}}(th_r) \Psi_n(t) \, \mathrm{d}t$. Here, for all $r$, $h_r > 0$ is a bandwidth parameter; in particular, $h_k$ needs to be chosen to ensure consistency of the estimator of $f_X$. We choose $h_r$ that minimizes the asymptotic mean squared error (AMSE) of the estimator $\widehat{\theta}_r$. As in the homoscedastic case, the AMSE can be decomposed as the sum of a squared bias term and a variance term, where, under sufficient conditions (see Delaigle and Meister (2007)), the latter is negligible; $h_r$ can thus be chosen on the basis of the sole asymptotic bias, given by

$$
\begin{aligned}
\mathrm{ABias}[\widehat{\theta}_r] = {} & (-1)^{k/2} \frac{2h_r^k}{k!} \mu_{K,k} \theta_{r+k/2} \\
& + \frac{1}{2\pi h_r^{2r+1}} \int t^{2r} |K^{\mathrm{ft}}(t)|^2 \bigg/ \left( \sum_{k=1}^n |f_{\varepsilon_k}^{\mathrm{ft}}(t/h_r)|^2 \right) \mathrm{d}t.
\end{aligned}
\tag{5.1}
$$

The procedure of Delaigle and Gijbels (2004) involves estimation of $\theta_{2k}$ by an estimator $\widehat{\theta}_{2k} = (4k)! / ((2\widehat{\sigma}_X)^{4k+1} (2k)! \pi^{1/2})$, obtained by assuming that $f_X$ is a normal density. Here, $\widehat{\sigma}_X$ is an estimator of the standard deviation of $X$, which, in our context, can be, for example, $\widehat{\sigma}_X^2 = [n^{-1} \sum_{i=1}^n Y_i^2 - (n^{-1} \sum_{i=1}^n Y_i)^2] - [n^{-1} \sum_{i=1}^n \mathrm{E}(\varepsilon_i^2) - (n^{-1} \sum_{i=1}^n \mathrm{E}(\varepsilon_i))^2]$.

## 5.2. Simulation results

We applied our estimator (2.2) to simulated examples from two densities $f_X$: (1) $X \sim 0.5N(-3, 1) + 0.5N(2, 1)$ and (2) $0.75N(0, 1) + 0.25N(1.5, 1/81)$. We considered four heteroscedastic models: (i) $\varepsilon_1, \ldots, \varepsilon_{n/2} \sim N(0, \sigma_1^2)$ and $\varepsilon_{n/2+1}, \ldots, \varepsilon_n \sim \mathrm{Laplace}(\sigma_2)$; (ii) $\varepsilon_1, \ldots, \varepsilon_{n/2} \sim N(0, \sigma_1^2)$ and $\varepsilon_{n/2+1}, \ldots, \varepsilon_n \equiv 0$; (iii) one error density $f_\varepsilon \sim N(0, \sigma_1^2)$, but a different number of replicated observations – here, we use the averaged data as in Section 3.1; and (iv) $\varepsilon_i \sim N(0, \sigma_3^2(1 + i/n))$. These are non-trivial situations because the target densities $f_X$ are not easy to estimate and normal errors are hard to deconvolve.

For density (1) (resp., density (2)), we took $\sigma_1$ and $\sigma_2$ such that $\mathrm{Var}(\varepsilon_i) = 25\%$ (resp., 10%) $\times \mathrm{Var}(X)$ and $\sigma_3^2 = 10\%$ (resp., 5%) $\times \mathrm{Var}(X)$. In each case, we generated 500 contaminated samples of size $n = 50$, 100 or 250 from the distribution of density (1) or (2). For each sample, we constructed the estimator (2.2) using the plug-in bandwidth of Section 5.1 and the kernel $K_2$. To evaluate performance, we calcu-lated, on a grid of 81 equidistant values of $x$, the quantiles $q_p(x)$ of the 500 estimates $\hat{f}_n(x)$ for $p = 0.1$, $0.25$, $0.5$, $0.75$ and $0.9$. In the graphs, we refer to $q_{0.5}$ as the median, $q_{0.25}$ and $q_{0.75}$ as the quartiles and $q_{0.1}$, $q_{0.9}$ as the deciles. We only present partial results, but our conclusions were also supported by the unreported cases.

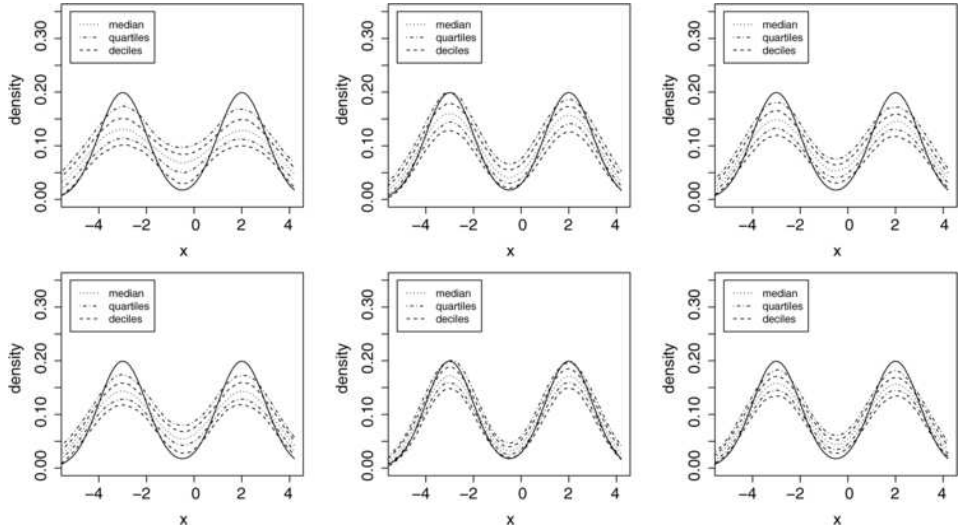

**Figure 1.** Estimators of density (1) from samples of size $n = 100$ (first row) and $n = 250$ (second row), generated from model (i) (left panel), (ii) (center panel) and (iv) (right panel).

In Figure 1, we show some quantile curves constructed from samples of size $n = 100$ and 250, generated from density (1) under models (i), (ii) and (iv). As expected by the theory, these graphs show a clear improvement of the results from (i) to (ii) and when the sample size increases. We also see that our method does not have particular problems in dealing with the case of individual errors.

In Figure 2, we compare the results for density (2) and samples of size $n = 100$ or 250 coming from models (i), (ii) and (iii), where 25% of the observations are not replicated and 50% (resp., 25%) of the observations are replicated twice (resp., ten times). Here, again, we see an improvement of the quality of the estimator from model (i) to model (ii) and the estimator handles the case of a different number of replicated measurements without any particular difficulty.

Additional results not reported here (see Delaigle and Meister (2007)) showed that the data-driven bandwidth procedure suffers from only a small loss of performance compared to the optimal bandwidth. In addition, although, asymptotically, the estimator that discards the observations contaminated by the smoothest errors has the same behavior as the estimator that uses all the observations, the latter had better practical properties, especially for the smallest sample sizes. Finally, our method worked considerably better than the one that ignores the errors in the data.

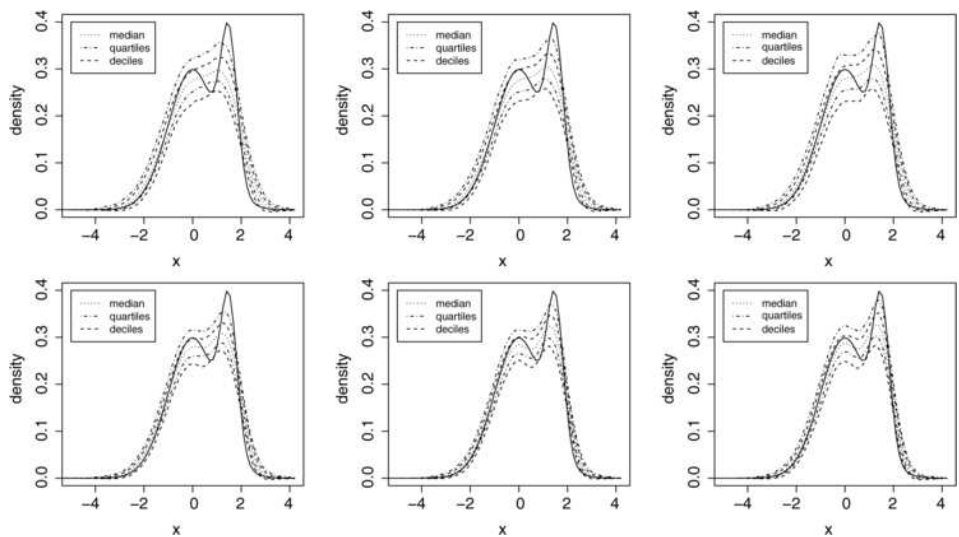

**Figure 2.** Estimators of density (2) from samples of size $n = 100$ (first row) and $n = 250$ (second row), in the case of normal and Laplace errors (first column), partially normally contaminated (second column) and replicated observations with normal errors (third column).

## 6. Proofs

**Proof of Theorem 2.1.** Part (a) follows from (C.2), (C.3), (C.5) and (2.5) applied to the fact that the MISE of the estimator is bounded by the sum of the first two terms of (2.3), which, in turn, is bounded by

$$\sup_{f_X \in \mathcal{F}_{\beta,C}} \mathrm{MISE}_n(f_X) = O\left(\int_0^{\omega_n} \left[\sum_{j=1}^n |f_{\varepsilon_j}^{\mathrm{ft}}(t)|^2\right]^{-1} \mathrm{d}t, \omega_n^{-2\beta}\right). \tag{6.1}$$

Concerning part (b), we note that Fan (1991a, b, 1993) derives theoretical lower bounds for standard density deconvolution under Hölder conditions; those results can be extended to Sobolev classes (see Neumann (1997)). Since we are considering a problem with non-identically distributed data, a new concept is required.

Let $f_0(x) = \pi^{-1}(1+x^2)^{-1}$ be the Cauchy density and set $f_1(x) = (1-\cos x)/(\pi x^2)$ with $f_1^{\mathrm{ft}}(t) = (1-|t|) \cdot 1_{[-1,1]}(t)$. We introduce the densities

$$f_\theta(x) = \tfrac{1}{2} f_0(x) + \tfrac{1}{2} f_1(x) + \sum_{j=\lfloor m_n \rfloor}^{2\lfloor m_n \rfloor} j^{-\beta-(1/2)} \theta_j \cos(2jx) f_1(x),$$

with $\theta_j \in \{0,1\}$. For $C$ and $n$ large enough, all $f_\theta$'s are contained in $\mathcal{F}_{\beta,C}$. Similarly to Fan ([1993](#)), we randomize the vector $\theta$ so that the $\theta_j$'s are i.i.d. with $P(\theta_j = 0) = 1/2$ and define $\theta_{j,0} = (\theta_{\lfloor m_n \rfloor}, \dots, \theta_{j-1}, 0, \theta_{j+1}, \dots, \theta_{2\lfloor m_n \rfloor})$ and $\theta_{j,1}$ accordingly. An application of Parseval's identity, combined with the fact that the $f_1^{\mathrm{ft}}(\cdot - 2j)$'s, $j$ integer, have disjoint supports, shows that after calculating the expectation with respect to $\theta_j$, we obtain, for any estimator $\widehat{f}_n$, that

$$\mathrm{E}_\theta \mathrm{E}_{f_\theta} \|\widehat{f}_n - f_\theta\|_{L_2(\mathbb{R})}^2 \geq (2\pi)^{-1} \sum_{j=\lfloor m_n \rfloor}^{2\lfloor m_n \rfloor} \mathrm{E}_\theta \mathrm{E}_{f_\theta} \int_{2j-1}^{2j+1} |\widehat{f}_n^{\mathrm{ft}}(t) - f_\theta^{\mathrm{ft}}(t)|^2 \, dt$$

$$\geq \mathrm{const.} \sum_{j=\lfloor m_n \rfloor}^{2\lfloor m_n \rfloor} \int_{2j-1}^{2j+1} |f_{\theta_{j,0}}^{\mathrm{ft}}(t) - f_{\theta_{j,1}}^{\mathrm{ft}}(t)|^2 \, dt \qquad (6.2)$$

$$\geq \mathrm{const.} \sum_{j=\lfloor m_n \rfloor}^{2\lfloor m_n \rfloor} j^{-2\beta-1} \geq \mathrm{const.} \cdot m_n^{-2\beta}$$

if, for any $|j| \in [\lfloor m_n \rfloor, 2\lfloor m_n \rfloor]$ and any $\theta_l \in \{0,1\}$ with $l \neq j$, we have

$$\int \cdots \int \min\left(\prod_{k=1}^n h_{k;\theta_{j,0}}(y_k), \prod_{k=1}^n h_{k;\theta_{j,1}}(y_k)\right) dy_1 \cdots dy_n \geq \mathrm{const.} > 0, \qquad (6.3)$$

with the densities $h_{k;\theta_{j,\bullet}} = f_{\theta_{j,\bullet}} * f_{\varepsilon_k}$. By applying LeCam's inequality (see, e.g., Devroye ([1987](#)), page 7) and the logarithmic function to both sides of ([6.3](#)), we see that ([6.3](#)) is satisfied if

$$\sum_{k=1}^n (1 - a_{j,k,n})/a_{j,k,n} = O(1) \qquad (6.4)$$

holds for all $|j| \in [\lfloor m_n \rfloor, 2\lfloor m_n \rfloor]$, where we write

$$a_{j,k,n} := \int [(f_{\theta_{j,0}} * f_{\varepsilon_k})(x)(f_{\theta_{j,1}} * f_{\varepsilon_k})(x)]^{1/2} \, dx.$$

Due to $f_{\theta_{j,\bullet}} \geq (1/2) f_0$, we see that $a_{j,k,n} \geq 1/2$ and, hence, ([6.4](#)) follows from

$$\sum_{k=1}^n \chi^2(h_{k;\theta_{j,0}}, h_{k;\theta_{j,1}}) = O(1), \qquad (6.5)$$

where $\chi^2(f, g) := \int (f - g)^2 / f \, \mathrm{d}x$ denotes the $\chi^2$-distance of densities. This generalizes the condition in Fan (1991a, b), $\chi^2(h_{1;\theta_{j,0}}, h_{1;\theta_{j,1}}) = O(1/n)$, to the case of heteroscedastic contamination. We notice that the left-hand side of (6.5) is bounded above by

$$O(m_n^{-2\beta-1}) \sum_{k=1}^n \int \frac{[\cos(2j\cdot)f_1 * f_{\varepsilon_k}]^2(x)}{[f_0 * f_{\varepsilon_k}](x)} \, \mathrm{d}x. \tag{6.6}$$

Unlike in the situation of i.i.d. data, the denominator in (6.6) still depends on $k$ and $n$. Condition (C.1) annuls this difficulty as we have

$$[f_0 * f_{\varepsilon_k}](x) \geq \pi^{-1} \int_{|y| \leq \alpha} [1 + (x - y)^2]^{-1} f_{\varepsilon_k}(y) \, \mathrm{d}y$$

$$\geq \pi^{-1} [1 + 2\alpha^2 + 2x^2]^{-1} \int_{|y| \leq \alpha} f_{\varepsilon_k}(y) \, \mathrm{d}y$$

$$\geq \text{const.} \cdot [1 + x^2]^{-1}.$$

Therefore, applying the Fourier representation of the Sobolev norm, term (6.6) is bounded above by

$$O(m_n^{-2\beta-1}) \sum_{k=1}^n \int (|f_1^{\mathrm{ft}}(t - 2j) f_{\varepsilon_k}^{\mathrm{ft}}(t)|^2 + |f_1^{\mathrm{ft}\prime}(t - 2j) f_{\varepsilon_k}^{\mathrm{ft}}(t)|^2 + |f_1^{\mathrm{ft}}(t - 2j) f_{\varepsilon_k}^{\mathrm{ft}\prime}(t)|^2) \, \mathrm{d}t$$

$$\leq O(m_n^{-2\beta-1}) \sum_{k=1}^n |\overline{\varphi}_{k,n}(m_n)|^2,$$

due to (C.3) and (C.4). Finally, (2.5) implies (6.5), which proves the theorem. $\qquad\square$

**Proof of Theorem 3.1.** We introduce the orthonormal $r_{j,n} \times r_{j,n}$ matrices $A_{j,n}$ which consist of $r_{j,n}^{-1/2} \cdot (1, \ldots, 1)$ as their first row. Setting $W_{j,\bullet} := A_{j,n} Y_{j,\bullet}$ with $Y_{j,\bullet} := (Y_{j,1}, \ldots, Y_{j,r_{j,n}})^t$, we notice that $W_{j,1} = r_{j,n}^{1/2} \overline{Y}_j$, while the other components of $W_{j,\bullet}$ are measurable in the $\sigma$-algebra generated by $\varepsilon_{j,1}, \ldots, \varepsilon_{j,r_{j,n}}$ since any row of $A_{j,n}$ except the first one sums to zero, due to the orthonormal structure of $A_{j,n}$. Concerning the density $f_{Y_{j,\bullet}}$ of $Y_{j,\bullet}$, we derive

$$f_{Y_{j,\bullet}}(y_{j,\bullet}) = \left(\frac{1}{\sqrt{2\pi}\sigma}\right)^{r_{j,n}} \int f_X(x) \exp(-\|y_{j,\bullet} - (x + \mu) \cdot (1, \ldots, 1)^t\|^2 / (2\sigma^2)) \, \mathrm{d}x$$

$$= \left(\frac{1}{\sqrt{2\pi}\sigma}\right)^{r_{j,n}} \int f_X(x) \exp(-\|A_{j,n} y_{j,\bullet} - r_{j,n}^{1/2}(x + \mu) \cdot (1, 0, \ldots, 0)^t\|^2 / (2\sigma^2)) \, \mathrm{d}x$$

$$= \left(\frac{1}{\sqrt{2\pi}\sigma}\right)^{r_{j,n}} \cdot \exp\left(-\frac{1}{2\sigma^2} \sum_{k=2}^{r_{k,n}} |w_{j,k}|^2\right)$$

$$\times \int f_X(x) \exp(-|w_{j,1} - r_{j,n}^{1/2}(x+\mu)|^2/(2\sigma^2)) \, \mathrm{d}x,$$

where $\|\cdot\|$ denotes the Euclidean norm and $w_{j,\bullet} = A_{j,n}y_{j,\bullet}$. Therefore, we see that the conditional distribution of $Y_{j,\bullet}$ given $W_{j,1}$ and, hence, the distribution of all available data,

$$dP(Y_{\bullet,\bullet} = y_{\bullet,\bullet} \mid W_{\bullet,1}) = \prod_{j=1}^{n} dP(Y_{j,\bullet} = y_{j,\bullet} \mid W_{j,1}),$$

do not depend on $f_X$. Thus we have shown sufficiency and the proof is complete. $\qquad \square$

**Proof of Theorem 3.2.** First, we assume condition (3.2) and take the sinc kernel $K_1$. In view of (3.2), for an arbitrarily small $\gamma' \in (0, \gamma)$, we can choose $\alpha' \in (0, \alpha)$ sufficiently small so that $\liminf_{n\to\infty} \#J_{n,\gamma',\alpha'}/\#J_{n,\gamma,\alpha} \geq 1$. This shows that, in (3.2), we can choose $\gamma = \gamma'$ and $\alpha = \alpha'$ with $\alpha'/2 + \gamma' - \delta < 0$. Setting $\omega_n = n^{\delta/(2\beta+1)}$ for $\delta \leq (\beta+1/2)\alpha'$ and $\omega_n = n^{\alpha'/2}$ otherwise, we learn from (6.1) that the bias term converges at algebraic rates. The variance has the upper bound

$$O(\omega_n) \cdot \left( \sum_{j=1}^{n} \exp(-\sigma_{j,n}^2 \omega_n^2) \right)^{-1}$$

$$\leq O(\omega_n) \cdot \left( \sum_{j \in J_{n,\gamma',\alpha'}} \exp(-\sigma_{j,n}^2 \omega_n^2) \right)^{-1}$$

$$\leq O(\omega_n) \cdot \left( \sum_{j \in J_{n,\gamma',\alpha'}} n^{-\gamma'} \right)^{-1}$$

$$\leq O(\omega_n n^{\gamma'-\delta}) \leq O(n^{\alpha'/2+\gamma'-\delta}).$$

Hence, the algebraic decay of the MISE has been established.

For the reverse implication, assume that the supremum of the MISE (and thus the bias and the variance terms in (6.1)) converges with an algebraic rate. The bias term then implies that $\omega_n \geq c \cdot n^s$ with $s > 0$, while the variance term is bounded below by

$$\mathrm{const.} \cdot \frac{n^s}{2} \cdot \left( \sum_{j=1}^{n} \exp(-\sigma_{j,n}^2 n^{2s}/4) \right)^{-1}$$

$$= \mathrm{const.} \cdot n^s \cdot \left( \sum_{j \in J_{n,4,2s}} \exp(-\sigma_{j,n}^2 n^{2s}/4) + \sum_{j \in J_{n,4,2s}^c} \exp(-\sigma_{j,n}^2 n^{2s}/4) \right)^{-1}$$

$$\geq \mathrm{const.} \cdot n^s \cdot (\#J_{n,4,2s} + n^{-1} \cdot \#J_{n,4,2s}^c)^{-1} \geq \mathrm{const.} \cdot n^s \cdot (\#J_{n,4,2s} + 1)^{-1}.$$

We deduce the existence of a $\delta > 0$ and a $c > 0$ such that $\#J_{n,4,2s} \geq c \cdot n^\delta$. $\qquad \square$

**Proof of Theorem 3.3.** (a) From (3.3), we can construct a sequence $(\omega_n)_n \to \infty$ such that

$$\omega_n \left( \sum_{j=1}^n \exp(-\sigma_{j,n}^\gamma \omega_n^\gamma) \right)^{-1} \overset{n \to \infty}{\longrightarrow} 0$$

for any known parameters $\sigma_{j,n}$. It follows that the variance term of estimator (2.2) converges to 0 due to Lemma 2.1 – as does the bias term as $\omega_n \to \infty$.

(b) We assume that (3.3) does not hold. Then, there exist $\omega_0 > 0$ and $M > 0$ such that

$$\sum_{j=1}^n \exp(-\sigma_{j,n}^\gamma \omega_0^\gamma) \leq M \tag{6.7}$$

for infinitely many $n$. In the sequel, we restrict our consideration to those $n$. We may assume $\omega_0$ to be arbitrarily large without affecting the validity of (6.7) and we also note that only a bounded number of the $\sigma_{j,n}$'s can be less than 1. Hence, in the view of the asymptotic behavior, we may assume that $\sigma_{j,n} > 1$, without loss of generality. For any $\omega_1 > \omega_0$, we have

$$\sum_{j=1}^n \exp(-\sigma_{j,n}^\gamma \omega_1^\gamma / 4) \leq M \exp(\omega_0^\gamma - \omega_1^\gamma / 4). \tag{6.8}$$

We introduce the density $f$ with Fourier transform $f^{\mathrm{ft}}(t) = (1 - |t/(2\omega_1)|) \cdot 1_{[-2\omega_1, 2\omega_1]}(t)$ and the density $\tilde{f}$ whose Fourier transform is supported on $[-3\omega_1, 3\omega_1]$ and coincides with $f^{\mathrm{ft}}(t)$ on its restriction to $[-\omega_1, \omega_1]$. On $[\omega_1, 3\omega_1]$, the even function $\tilde{f}^{\mathrm{ft}}(t)$ is defined as the linear connection of the points $(3\omega_1, 0)$ and $(\omega_1, f^{\mathrm{ft}}(\omega_1))$. The existence of $\tilde{f}$ is guaranteed by Pólya's criterion (see Lukacs (1970), page 83, Theorem 4.3.1). We notice that $f, \tilde{f} \in \mathcal{F}_{\beta, C}$ for any $\beta > 1/2$, with $C$ sufficiently large. The Parseval identity gives us

$$\|f - \tilde{f}\|_{L_2(\mathbb{R})}^2 \geq \omega_1 / (48\pi).$$

Equipped with those results, we fix an arbitrary estimator $\widehat{f}_n$ of $f_X$ and consider

$$\begin{aligned}
&\mathrm{E}_f \|\widehat{f}_n - f\|_{L_2(\mathbb{R})}^2 + \mathrm{E}_{\tilde{f}} \|\widehat{f}_n - \tilde{f}\|_{L_2(\mathbb{R})}^2 \\
&\quad \geq \mathrm{E}_f \|\widehat{f}_n - f\|_{L_2(\mathbb{R})}^2 + \mathrm{E}_f \|\widehat{f}_n - \tilde{f}\|_{L_2(\mathbb{R})}^2 \\
&\qquad - |\mathrm{E}_f \|\widehat{f}_n - \tilde{f}\|_{L_2(\mathbb{R})}^2 - \mathrm{E}_{\tilde{f}} \|\widehat{f}_n - \tilde{f}\|_{L_2(\mathbb{R})}^2| \\
&\quad \geq \omega_1 / (48\pi) - O\left( \sum_{j=1}^n \|(f - \tilde{f}) * f_{\varepsilon_j}\|_{L_1(\mathbb{R})} \right).
\end{aligned} \tag{6.9}$$

Therefore, we can establish inconsistency by showing that (6.9) is bounded away from zero for a fixed choice of $\omega_1 > 0$. To this end, we need an upper bound for each $\|(f -$ $\tilde f) * f_{\varepsilon_j}\|_{L_1(\mathbb{R})}$. Employing the Cauchy density $f_0(x) = [\pi(1+x^2)]^{-1}$, we use the Cauchy–Schwarz inequality to obtain

$$\|(f - \tilde f) * f_{\varepsilon_j}\|_{L_1(\mathbb{R})} \leq \left(\pi \int |[(f - \tilde f) * f_{\varepsilon_j}](x)|^2 (1+x^2)\,\mathrm{d}x\right)^{1/2}. \qquad (6.10)$$

As in the proof of Theorem 2.1, the Fourier representation of the Sobolev norm leads to the following upper bound for the right-hand side of (6.10):

$$\left(\int [|(f^{\mathrm{ft}}(t) - \tilde f^{\mathrm{ft}}(t)) f_{\varepsilon_j}^{\mathrm{ft}}(t)|^2 + |(f^{\mathrm{ft}\prime}(t) - \tilde f^{\mathrm{ft}\prime}(t)) f_{\varepsilon_j}^{\mathrm{ft}}(t)|^2 + |(f^{\mathrm{ft}}(t) - \tilde f^{\mathrm{ft}}(t)) f_{\varepsilon_j}^{\mathrm{ft}\prime}(t)|^2\,\mathrm{d}t\right)^{1/2}.$$

Therefore, we see that (6.9) has the lower bound

$$\omega_1/(48\pi) - O\left(\sum_{j=1}^n \exp(-\sigma_{j,n}^\gamma \omega_1^\gamma/4)\right) \qquad (6.11)$$

when selecting $\omega_1$ sufficiently large. We apply (6.8) so that for appropriate constants $c_1, c_2 > 0$, (6.11) is bounded below by $c_1\omega_1 - c_2\exp(-\omega_1^\gamma/4)$. Choosing $\omega_1 > 0$ large enough, while $\omega_0$ is fixed, guarantees a positive lower bound for (6.9) and, hence, inconsistency. □

# Acknowledgements

A. Delaigle's research was supported by a Hellman Fellowship and a Maurice Belz Fellowship. A. Meister's research was partly carried out at the Centre for Mathematics and its Applications, Australian National University, Canberra, Australia.